\documentclass{amsart}
\usepackage[all]{xy}
\usepackage{epsfig,cite}
\usepackage{amssymb,amsmath, amsthm}
\usepackage{times}

\begin{document}
\sloppy \raggedbottom
\setcounter{page}{1}

\newpage
\setcounter{figure}{0}
\setcounter{equation}{0}
\setcounter{footnote}{0}
\setcounter{table}{0}
\setcounter{section}{0}
\newtheorem{thm}{Theorem}[subsection]
\newtheorem{cor}[thm]{Corollary}
\newtheorem{thm-def}[thm]{Theorem-Definition}
\newtheorem{lem}[thm]{Lemma}
\newtheorem{prop}[thm]{Proposition}
\newtheorem{defn}[thm]{Definition}
\theoremstyle{definition}
\newtheorem{rem}[thm]{Remark}
\numberwithin{equation}{section}

\newcommand{\Spec}{{\rm Spec}}

\newcommand{\tensor}{\otimes}

\newcommand{\vacA}{\mathbf{1}}
\newcommand{\vacB}{|0\rangle}
\newcommand{\vac}{\vacB}

\newcommand{\Real}{\mathbb R}
\newcommand{\RR}{\mathbb R}
\newcommand{\Cplx}{\mathbb C}
\newcommand{\CC}{\mathbb C}
\newcommand{\ZZ}{\mathbb Z}
\newcommand{\LL}{\mathbb L}
\newcommand{\VV}{\mathbb V}
\newcommand{\MM}{\mathbb M}
\newcommand{\PP}{\mathbb P}
\newcommand{\WW}{\mathbb W}
\newcommand{\vacu}{\mathbf 1}

\newcommand{\Field}{\mathbb F}
\newcommand{\Rat}{\mathbb Q}

\newcommand{\gr}{{\rm gr\ }}
\newcommand{\grG}{{\rm gr^{\mathcal{G}}\ }}
\newcommand{\grF}{{\rm gr^{\mathcal{F}}\ }}
\newcommand{\grH}{{\rm gr^{\mathcal{H}}\ }}

\newcommand{\pole}{\Cplx}

\newcommand{\de}{\partial}
\newcommand{\oX}{{\mathcal{O}}_{X}}

\newcommand{\calO}{{\mathcal{O}}}

\newcommand{\calA}{{\mathcal{A}}}
\newcommand{\calB}{{\mathcal{B}}}
\newcommand{\calC}{{\mathcal{C}}}
\newcommand{\calD}{{\mathcal{D}}}
\newcommand{\calE}{{\mathcal{E}}}

\newcommand{\calF}{{\mathcal{F}}}
\newcommand{\calG}{{\mathcal{G}}}
\newcommand{\calH}{{\mathcal{H}}}

\newcommand{\calK}{{\mathcal{K}}}

\newcommand{\calT}{{\mathcal{T}}}

\newcommand{\cA}{\mathcal{A}}
\newcommand{\cB}{\mathcal{B}}

\newcommand{\cD}{{\mathcal{D}}}

\newcommand{\cF}{{\mathcal{F}}}
\newcommand{\cG}{{\mathcal{G}}}
\newcommand{\cH}{{\mathcal{H}}}
\newcommand{\cJ}{{\mathcal{J}}}
\newcommand{\cI}{{\mathcal{I}}}
\newcommand{\cL}{{\mathcal{L}}}
\newcommand{\cM}{{\mathcal{M}}}
\newcommand{\cN}{\mathcal{N}}
\newcommand{\cO}{{\mathcal{O}}}
\newcommand{\cR}{{\mathcal{R}}}
\newcommand{\cS}{{\mathcal{S}}}

\newcommand{\cT}{{\mathcal{T}}}

\newcommand{\cV}{{\mathcal{V}}}
\newcommand{\cW}{{\mathcal{W}}}

\newcommand{\cZ}{{\mathcal{Z}}}
\newcommand{\cE}{{\mathcal{E}}}
\newcommand{\cC}{{\mathcal{C}}}

\newcommand{\TT}{\partial} 

\newcommand{\wt}[1]{\Delta_{#1}}

\newcommand{\op}[1]{{#1}_{(-1)} }

\newcommand{\opp}[2]{ {#1}_{(  #2  )} }

\newcommand{\zero}{{}_{(0)} }
\newcommand{\one}{{}_{(1)} }
\newcommand{\opm}{{}_{(-1)} }
\newcommand{\ops}[1]{{}_{({#1})} }


\newcommand{\Hom}{{\rm Hom}}
\newcommand{\Aut}{{\rm Aut}}
\newcommand{\Endom}{{\rm End\,}}
\newcommand{\End}{{\mathrm End\,}}

\newcommand{\Iso}{{\rm Iso}}

\newcommand{\id}{{\rm id\ }}

\newcommand{\Ind}{{\rm Ind\ }}
\newcommand{\im}{{\rm Im\ }}

\newcommand{\spn}{{\rm span}}
\newcommand{\Der}{{\rm Der\ }}

\newcommand{\iso}{\stackrel{\sim}{\rightarrow}}

\newcommand{\res}{{\rm res}}
\newcommand{\Resz}{{\rm Res}_z}
\newcommand{\ta}{{\tilde{a}}}
\newcommand{\tb}{{\tilde{b}}}
\newcommand{\ttt}{{\tilde{t}}}
\newcommand{\ts}{{\tilde{s}}}
\newcommand{\tg}{{\tilde{g}}}
\newcommand{\tf}{{\tilde{f}}}

\newcommand{\txi}{{\tilde{\xi}}}
 \newcommand{\teta}{{\tilde{\eta}}}

\newcommand{\tlam}{{\tilde{\lambda}}}

\newcommand{\la}{\lambda}

\newcommand{\fg}{\frak{g}}
\newcommand{\fa}{\frak{a}}
\newcommand{\fb}{\frak{b}}
\newcommand{\fn}{\frak{n}}
\newcommand{\fh}{\frak{h}}
\newcommand{\fz}{\frak{z}}
\newcommand{\fzl}{\fz^{\text{reg},\lambda}}

\newcommand{\fU}{{\frak{U}}}
\newcommand{\fM}{{\frak{M}}}

\newcommand{\sch}[2]{ \Gamma ({#1}, {#2} ) }

\newcommand{\catleq}{{\cV}ert^{\leq 1}}
\newcommand{\slt}{\frak{s}\frak{l}_2}
\newcommand{\affsl}{\widehat{\frak{s}\frak{l}}_2}

\newcommand{\prline}{\mathbb{P}^1}

\newcommand{\tdl}{TD^{\lambda}}
\newcommand{\Vcrit}{V_{\kappa_c}(\slt)}
\newcommand{\Symb}{{\rm Symb\  }}

\newcommand{\omcl}{\Omega_X^{1,\,cl}}
\newcommand{\homcl}{H^1 (X, \omcl)}
\newcommand{\funomega}{Fun(\homcl)}


\newcommand{\cdo}{{\cD}_X^{ch}}
\newcommand{\dl}{\cD^\cL}
\newcommand{\twcdo}{{\cD}_X^{ch, tw}}
\newcommand{\twcdoloc}{{\stackrel{\circ}{\cD}}_X^{ch, tw}}
\newcommand{\gtdo}{{\cD}_X^{tw}}
\newcommand{\cdofun}[1]{\cD_{#1}\otimes\funomega}

\newcommand{\Pic}{ {\rm Pic}  }
\newcommand{\cch}{\Cplx_{\chi(z)}}
\newcommand{\cps}{\Cplx_{\Psi}}
\newcommand{\dchmod}{\cdo{\textrm{-}}Mod }
\newcommand{\dchtwmod}{\twcdo{\textrm{-}}Mod }
\newcommand{\dmod}{\cD_X{\textrm{-}}Mod}
\newcommand{\dlmod}{\cD^{\cL}_X{\textrm{-}}Mod}
\newcommand{\dlm}{\cD^{\cL}_X{\textrm{-}}Mod}

\newcommand{\Dchtwmodchi}{\twcdo{\textrm{-}}Mod_\chi }

\newcommand{\dlamod}{\gtdo{\textrm{-}}Mod}
\newcommand{\tensL}{\otimes_{\oX}\!\cL}
\newcommand{\ttens}{\overline{\otimes}}
\newcommand{\bM}{\bar{M}}
\newcommand{\bg}{\bar{g}}
\newcommand{\btau}{\bar{\tau}}
\newcommand{\bxi}{\bar{\xi}}
\newcommand{\ba}{\bar{a}}

\newcommand{\abs}[1]{\left\vert#1\right\vert}
\newcommand{\set}[1]{\left\{#1\right\}}
\newcommand{\eps}{\varepsilon}
\newcommand{\To}{\longrightarrow}
\newcommand{\BX}{\mathbf{B}(X)}


\newcommand{\CDOP}{\cD^{ch}_{\prline}}
\newcommand{\TWCDOP}{\cD^{ch,tw}_{\prline}}
\newcommand{\twcdop}{\cD^{ch,tw}_{\prline}}

\newcommand{\basis}{ \set{\tau^i} }
\newcommand{\basisom}{\{ \omega_i\}}
\newcommand{\Lie}{{\rm Lie}}

\newcommand{\cfdeg}{{\rm conf.deg}}

\newcommand{\tP}{\tilde{P}}
\newcommand{\tQ}{\tilde{Q}}

\newcommand{\partialtau}[1]{D_{\de_{#1}}}

\newcommand{\partialom}[1]{D_{\omega_{#1}}}

\newcommand{\Fields}{{Fields}}
\newcommand{\ghat}{{\hat{\fg}}}
\newcommand{\qf}[2]{\langle {#1}, {#2}  \rangle}
\newcommand{\vlk}{v_{\lambda, k}}
\newcommand{\uH}{\underline{\rm{H}}}

\pagestyle{plain}


\title{A Vertex Algebra Attached to the Flag Manifold and Lie Algebra Cohomology}



\author{T.Arakawa \and F.Malikov}
\thanks{T.A.\ is partially  supported
by the JSPS Grant-in-Aid  for Scientific Research (B)
No.\ 20340007.}
\thanks{F.M. is partially supported by an NSF grant.}

\maketitle



\begin{abstract}
We compute the cohomology vertex algebra $H^\bullet(G/B,\cD^{ch}_{G/B})$ to the effect
that it is isomorphic to $H^\bullet(\fn,\CC)\otimes \VV(\fg)_{\chi=0}$. We then find the
Friedan-Martinec-Shenker-Borisov  bosonization of $H^\bullet(\prline,\cD^{ch}_{\prline})$
and verify that the latter algebra vanishes  nonperturbatively.
\end{abstract}


\section{Introduction}
\label{introduction}

If $X$ is a smooth algebraic variety and $\cV$ a sheaf of vertex algebras over it, then
$H^\bullet(X,\cV)$ is naturally a vertex algebra. Although various examples of such sheaves,
most notably sheaves of algebras of chiral differential operators (CDO), \cite{BD,MSV,GMSII} and references therein,
have been studied, almost no
cohomology vertex algebras have been explicitly described, except for the case where $X=\PP^n$, \cite{MS}.

In the present paper we fill this gap by determining the vertex algebra structure of $H^\bullet(G/B,\cD^{ch}_{G/B})$,
where $G/B$ is a flag manifold and $\cD^{ch}_{G/B}$ is the CDO over it, \cite{MSV,GMSII}. The CDO $\cD^{ch}_{G/B}$ is a sheaf
of $\widehat{\fg}$-modules at the critical level, and the $\widehat{\fg}$-module structure of $H^\bullet(G/B,\cD^{ch}_{G/B})$
was computed in \cite{AM}. It is clear from \cite{FG3} that $H^\bullet(G/B,\cD^{ch}_{G/B})$ is a direct sum of several copies of  $\VV(\fg)_{\chi=0}$, the vacuum $\widehat{\fg}$-module at the critical level and zero central character. The computation
carried out in \cite{AM} establishes a natural 1-1 correspondence between these copies
and Schubert cells of $G/B$. This suggests that as a vertex algebra $H^\bullet(G/B,\cD^{ch}_{G/B})$ is isomorphic to $H^\bullet(G/B,\CC)\otimes \VV(\fg)_{\chi=0}$, a proposition  immediately seen to be wrong because it distorts the cohomological degree. What happens in reality is rather a manifestation of Bott's ``strange equality'', \cite{Bott}:
\[\text{dim}H^i(\fn,\CC)=\text{dim}H^{2i}(G/B,\CC),
\]
where $\fn\subset\fg$ is a maximal nilpotent subalgebra. We show that the associative commutative algebra $H^\bullet(\fn,\CC)$ carries a derivation -- thus becoming a commutative vertex algebra -- so that $H^\bullet(G/B,\cD^{ch}_{G/B})$ is isomorphic to $H^\bullet(\fn,\CC)\otimes \VV(\fg)_{\chi=0}$ as a vertex algebra. Furthermore,
$H^\bullet(\fn,\CC)$ thus embedded in $H^\bullet(G/B,\cD^{ch}_{G/B})$ is the center of $H^\bullet(G/B,\cD^{ch}_{G/B})$.

Our result is not without omission for we do not compute the derivation, see Theorem~\ref{main-thm}. We show, however, that this derivation is 0 when $\fg$ is either $sl_2$ or $sl_3$ and expect it to be always 0.

Next, we focus on the case of $sl_2$ and establish what can be called a {\it
Friedan-Martinec-Shenker-Borisov  bosonization} of $H^\bullet(\prline,\cD^{ch}_{\prline})$. This means a differential graded vertex algebra $(V_L^0,D)$ such that the cohomology vertex algebra $H_D(V_L^0)$ is $H^\bullet(\prline,\cD^{ch}_{\prline})$. The DGVA $(V_L^0,D)$ can be included in a flat family $(V_L^t,D)$, $t\in\CC$, where $V_L^t$ is a certain lattice vertex algebra if $t\neq0$. We show that generically the cohomology vanishes:
$H_D(V_L^t)=0$ if $t\neq0$. This result should be contrasted with an analogous result of
\cite{MS}, where a similar computation gives that $H_D(V_Q^0)$ is a full-fledged infinite
dimensional vertex algebra, the cohomology of the chiral de Rham complex over a projective
space, while generically, if $t\neq0$, $H_D(V_Q^t)$ is the quantum cohomology of  the same projective space. If we are to believe that the passage from $(V_L^0,D)$ to $(V_L^t,D)$
is a way to take into account the instanton corrections (arguments in \cite{FL} seem to be in favor of this), then our result is a confirmation of Witten's prediction \cite{W} that
nonperturbatively the algebra $H^\bullet(X,\cD^{ch}_X)$ vanishes\footnote{a string theory verification of Witten's observation for any
flag manifold was obtained in \cite{MY1}}. 

\section{Cohomology vertex algebra}
\label{Cohomology vertex algebra}
\subsection{Multiplicative structures on cohomology}
Let $X$ be a smooth algebraic variety over $\CC$, $\cA$ a sheaf of $\CC$-vector
spaces over $X$. We would like to show that any multiplication on $\cA$, i.e.,
a sheaf morphism
\begin{equation}
\label{mult-on-sheaf}
m:\cA\otimes_{\CC}\cA\rightarrow\cA
\end{equation}
induces a canonical multiplication on cohomology
\begin{equation}
\label{mult-on-sheaf-coho}
h(m):H^{\bullet}(X,\cA)\otimes_{\CC}H^{\bullet}(X,\cA)\rightarrow H^{\bullet}(X,\cA).
\end{equation}
This and the discussion to follow should be termed well-known, but we failed to find a reference with
the necessary results. The earliest paper we were able to locate where such issues are analyzed is \cite{GH}.
Our thinking about the subject was greatly influenced by the Beilinson-Drinfeld approach to chiral algebras \cite{BD}; more on this in
Remark~\ref{pseudo-ten-rem}.

Consider the diagonal $\Delta: X\hookrightarrow X\times X$, the 2 projections
$p_i:X\times X\rightarrow X$, $i=1,2$, and the exterior tensor square
$\cA\boxtimes\cA=p_1^\bullet(\cA)\otimes_{\CC}p_2^\bullet(\cA)$, where $\Delta^\bullet$ is the
inverse image in the category of sheaves of $\CC$-vector spaces.
Since $\cA\otimes_{\CC}\cA=\Delta^{\bullet}(\cA\boxtimes\cA)$ and $\Delta^\bullet$ is the right adjoint
of $\Delta_\bullet$, the
direct image in the category of sheaves of $\CC$-vector spaces, we obtain
\[
\text{Hom}(\cA\otimes_{\CC}\cA,\cA) \stackrel{\sim}{\rightarrow}\text{Hom}(\cA\boxtimes\cA,\Delta_\bullet(\cA)).
\]
Hence (\ref{mult-on-sheaf})  gives
\begin{equation}
\label{mult-on-sheaf-adj}
\tilde{m}:\cA\boxtimes\cA\rightarrow\Delta_\bullet(\cA).
\end{equation}
Let us compute the cohomology of the sheaves involved. If $\cI$ is an injective module over $X$,
then $\Delta_\bullet(\cI)$ is an injective module over $X\times X$ with $\Gamma(X\times X,\Delta_\bullet(\cI))
=\Gamma(X,\cI)$. Therefore, if $\cA\rightarrow\cI^\bullet$
is an injective resolution, then $\Delta_\bullet(\cA)\rightarrow\Delta_\bullet(\cI^\bullet)$ is also, and we obtain
\begin{equation}\label{coho-diag}
H^\bullet(X\times X,\Delta_\bullet(\cA))=H^\bullet(X,\cA).
\end{equation}

Even if $\cI$ and $\cJ$ are injective, $\cI\boxtimes\cJ$ does not have to be, but its higher cohomology
vanish. Indeed, denoting by $\pi$ the projection to a point and factoring out $\pi= \pi\circ p_1$, we
obtain
\begin{eqnarray}
R  \pi(\cI\boxtimes\cJ)&=&\label{inj-ext-prod-coho}\\
R  \pi\circ Rp_{1\bullet}(\cI\boxtimes\cJ)&=&R\pi(\cI)\otimes_{\CC}\Gamma(X,\cJ)=
\Gamma(X,\cI)\otimes_{\CC}\Gamma(X,\cJ)
\nonumber.
\end{eqnarray}
It follows that the resolution $\cA\boxtimes\cA\rightarrow\cI^\bullet\boxtimes\cI^\bullet$ is adapted to the
computation of the cohomology $H^\bullet(X\times X,\cA\boxtimes\cA)$; (\ref{inj-ext-prod-coho}) implies
\begin{equation}
\label{coho-of-prod}
H^\bullet(X\times X,\cA\boxtimes\cA)=\Gamma(X,\cI^\bullet\boxtimes\cI^\bullet)=H^\bullet(X,\cA)\otimes_{\CC}H^\bullet(X,\cA).
\end{equation}
Multiplication (\ref{mult-on-sheaf-adj}) gives a multiplication on the resolution, well-defined
in the derived category,
\begin{equation}
\label{mult-on-resol}
\cI^\bullet\boxtimes\cI^\bullet\rightarrow\Delta_\bullet(\cI^\bullet).
\end{equation}
Taking the spaces of global
sections of both sides of (\ref{mult-on-resol}) and then using (\ref{coho-diag}) and
(\ref{coho-of-prod}), we obtain the desired multiplication
\begin{equation}
\label{def-of-mult}
h(m): H^\bullet(X,\cA)\otimes_{\CC}H^\bullet(X,\cA)\rightarrow H^\bullet(X,\cA).
\end{equation}

\subsection{The \u Cech complex and cup-product}
\label{The  ech complex and cup-product}

Suppose now  each $x\in X$ has a neighborhood $U\subset X$ so that $H^i(U,\cA|_{U})=0$ if $i> 0$.
Then for each fine enough open cover $\fU$, the (sheaf version of) \u Cech complex $\cC^\bullet=\cC^\bullet(\fU,\cA)$ is a resolution
of $\cA$. Hence a resolution
\begin{equation}
\cA\boxtimes\cA\rightarrow\cC^\bullet\boxtimes\cC^\bullet.
\end{equation}
Multiplication (\ref{mult-on-sheaf-adj}) allows us to define a {\em cup-product} on \u Cech cochains in a familiar manner
\begin{eqnarray}
\cC^\bullet\boxtimes\cC^\bullet&\rightarrow&\Delta_\bullet(\cC^\bullet),
\label{mult-cech-coch}\\
\cC^p\boxtimes\cC^q\ni f\boxtimes g&\mapsto& f\ast g,\;(f\ast g)(i_0,...,i_{p+q})=m(f(i_0,...,i_{p}),g(i_p,...,i_{q})).
\nonumber
\end{eqnarray}
We have obtained the following commutative diagram
\begin{displaymath}
\xymatrix{ \cA\boxtimes\cA\ar[d]\ar[r]&\cC^\bullet\boxtimes\cC^\bullet\ar[d]\\
 \Delta_\bullet(A)\ar[r]& \Delta_\bullet(\cC^\bullet)}
\end{displaymath}
Taking cohomology we obtain yet another, \u Cech, multiplication
\begin{equation}
\label{def-of-mult-cech}
\check{h}(m): H^\bullet(X,\cA)\otimes_{\CC}H^\bullet(X,\cA)\rightarrow H^\bullet(X,\cA).
\end{equation}
We wish to show that this multiplication coincides with the canonical one (\ref{def-of-mult}).
In fact, this is true at the cochain level. Note that $\cI^\bullet$ being injective we obtain
quasi-isomorphisms, $\cC^\bullet\rightarrow\cI^\bullet$, $\Delta_\bullet(\cC^\bullet) \rightarrow \Delta_\bullet(\cI^\bullet)$,
$\cC^\bullet\boxtimes\cC^\bullet\rightarrow\cI^\bullet\boxtimes\cI^\bullet$, canonical at the level of the derived category.

\begin{lem} Multiplications (\ref{mult-on-resol}) and (\ref{mult-cech-coch}) coincide in the derived category.
In particular,
\begin{equation}
\check{h}(m)=h(m).
\end{equation}
\end{lem}

{\em Proof.} Collect all the morphisms in sight in the following diagram:

\begin{displaymath}
\xymatrix{ \cA\boxtimes\cA\ar[ddd]\ar[rd]\ar[rr] & &\cI^\bullet\boxtimes\cI^\bullet\ar[ddd]\\
                                   &\cC^\bullet\boxtimes\cC^\bullet \ar[ur]\ar[d]\\
                                   &\Delta_\bullet(\cC^\bullet)\ar[rd]\\
 \Delta_\bullet(A)\ar[rr] \ar[ur]                      &  &\Delta_\bullet(\cI^\bullet)}
\end{displaymath}
By construction, commutative in this diagram are the top and bottom triangles, the left hand  trapezoid, and
the square. The right hand  trapezoid commutes only at the level of derived categories -- as indeed asserted by the lemma. To see this, one only needs
to precompose the 2 of its paths leading from $\cC^\bullet\boxtimes\cC^\bullet $ to $\Delta_\bullet(\cI^\bullet)$
and then use the listed commutativity properties. Applying the cohomology functor to this trapezoid, we obtain the commutative diagram
\begin{displaymath}
\xymatrix{ H^\bullet(X,\cA)\otimes_{\CC}H^\bullet(X,\cA)\ar@{=}[r]\ar[d]_{\check{h}(m)}&H^\bullet(X,\cA)\otimes_{\CC}H^\bullet(X,\cA)\ar[d]^{h(m)}\\
           H^\bullet(X,\cA)\ar@{=}[r]&H^\bullet(X,\cA)}
\end{displaymath}
The proof is complete. $\Box$

\subsection{Identities}
\label{Identities}
Assume that the algebra $\cA$ satisfies an order $n$ multilinearized identity. This means
there is a function $f: S_n\rightarrow\CC$, $S_n$ the symmetric group, such that
for each $U\subset X$ and each $n$-tuple $a_1,...,a_n\subset\cA(U)$
\begin{equation}
\label{identity-general}
\sum_{\sigma\in S_n,\,(( ))}f(\sigma)((x_{\sigma_1}\ast x_{\sigma_2}\ast\cdots x_{\sigma_n}))=0.
\end{equation}
We are not assuming associativity, and so $((x_{\sigma_1}\ast x_{\sigma_2}\ast\cdots x_{\sigma_n}))$ stands for the product of
the $x$'s taken in the indicated order and with the brackets inserted in some way that is somewhat schematically denoted by $((\,))$.

Each operation $\cA^{\otimes n}\rightarrow\cA$, $x_1\otimes\cdots \otimes x_n\mapsto ((x_1\ast x_2\cdots \ast x_n))$
extends to an operation on the resolution
$(\cI^\bullet)^{\boxtimes n}\rightarrow\Delta_\bullet^{(n)}(\cI^\bullet)$, where $\Delta^{(n)}:X\rightarrow X^{\times n}$ is the diagonal.

Each permutation $\cA^{\otimes n}\rightarrow\cA^{\otimes n}$ also extends to a permutation on $(\cI^\bullet)^{\boxtimes n}$, but in
order to obtain a morphism of complexes a
standard rule
of sign must be enforced.

If a map such as one defined by (\ref{identity-general}) equals 0, then the corresponding map on
$(\cI^{\bullet})^{\boxtimes n}$ is homotopic to 0 -- because a map such as (\ref{mult-on-resol}) is unique
up to homotopy. The passage to the cohomology proves the following.
\begin{lem}
\label{from-ident-sheaf-to-ident-coho}
If $\cA$ satisfies identity (\ref{identity-general}), then $H^\bullet(X,\cA)$ satisfies its superversion
\begin{equation}
\label{identity-general-coho}
f^{super}:\; \sum_{\sigma\in S_n}(-1)^{\sigma(x)}f(\sigma)((x_{\sigma_1}\ast x_{\sigma_2}\ast\cdots x_{\sigma_n}))=0,
\end{equation}
where $\sigma(x)$ is the number of inversions of odd cohomological degree $x$'s in $x_{\sigma_1}\ast x_{\sigma_2}\ast\cdots x_{\sigma_n}$.
\end{lem}

\begin{rem}
 \label{pseudo-ten-rem}
 All of this can be rephrased as follows.The  category with objects sheaves of vector spaces over $X$, morphisms various $Hom(\cA_1\boxtimes\cdots\boxtimes\cA_n,\Delta^{(n)}_{\bullet}(\cA))$ is a pseudo-tensor category \cite{BD}. An algebra $\cA$ defines
an algebra object in this category, the multiplication being a morphism
$\cA\boxtimes\cA\rightarrow\Delta_\bullet(\cA)$. If $\cA$ satisfies an order $n$ multilinearized identity, so will the corresponding object in the
pseudo-tensor category. The cohomology functor can be applied to
the multiplication $\cA\boxtimes\cA\rightarrow\Delta_\bullet(\cA)$ so as to endow $H^\bullet(X,\cA)$ with a $\CC$-algebra structure satisfying a super-version of the identity.
\end{rem}

\subsection{Vertex algebras}
\label{vertex-algebras}
A vertex algebra is a collection $(V,\vacu,_{(n)};n\in\ZZ)$, where
$V$ is a $\CC$-supervector space, $_{(n)}$ is a multiplication
\[
_{(n)}:V\otimes V\rightarrow V, a\otimes b\mapsto a_{(n)}b
\]
$\vacu\in V$ is a distinguished element known as the vacuum vector; the following axioms hold:
\begin{equation}
\label{nilpot}
a_{(n)}b=0\text{ if }n\gg 0,
\end{equation}
\begin{equation}
\label{vacuum}
\vacu_{(n)}=\delta_{n,-1}\text{Id},\; a_{(-1)}\vacu=a,\;a_{(n)}\vacu=0\text{ if }n\geq 0
\end{equation}

\begin{eqnarray}
   \label{Borcherds-identity}
   &\sum\limits_{j\geq 0} {m \choose j} (a\ops{n+j} b )\ops{m+k-j}c  = \\
   &\sum\limits_{j\geq 0} (-1)^{j} {n \choose j}\{ a\ops{m+n-j} (b\ops{k+j}c) - (-1)^{n}(-1)^{\tilde{a}\tilde{b}}b\ops{n+k-j} (a\ops{m+j}c) \}\nonumber
   \end{eqnarray}
If $\cV$ is a sheaf of vertex algebras, then each multiplication  is extended to an injective resolution of $\cV$
via  (\ref{mult-on-resol}) or to a \u Cech resolution via the cup-product (\ref{mult-cech-coch}). These multiplications define ones on the cohomology
$H^{\bullet}(X,\cV)$ via (\ref{def-of-mult}) or (\ref{def-of-mult-cech}).
The validity of the nilpotency condition (\ref{nilpot}) follows at once from the
cup-product formula (\ref{mult-cech-coch}).
The  multiplications on the cochains satisfy the Borcherds identity (\ref{Borcherds-identity}), an example of a trilinearized identity, up to homotopy, the  ones on the cohomology satisfy (\ref{Borcherds-identity}) on the nose; this
follows from sect.~\ref{Identities}, especially Lemma~\ref{from-ident-sheaf-to-ident-coho}, or rather its  obvious
multioperational version. We have established the following:
\begin{lem}
\label{vert-alg-str-on-coho}
If $\cV$ is a sheaf of vertex algebras, then $H^\bullet(X,\cV)$ carries
a canonical vertex algebra structure with  $\vacu\in\Gamma(X,\cV)$ the vacuum vector
and multiplications $_{(n)}$, $n\in\ZZ$, defined via (\ref{def-of-mult}) or (\ref{def-of-mult-cech}).
\end{lem}
Here is our main application of the above discussion.

\subsection{A CDO over the flag manifold}
\label{The cohomology vertex algebra of a CDO over the flag manifold}
Let $G$ be a simple complex Lie group, $B\subset G$ the Borel subgroup.  The flag
manifold $G/B$ carries an algebra of {\em chiral differential operators} (CDO),
an example of a sheaf of vertex algebras introduced initially in \cite{MSV}; see
also \cite{GMSII} for a different point of view.

The vertex algebra $\VV(\fg)_{-h^\vee}$ attached to the Lie algebra $\fg=\text{Lie}\,G$
at the critical level contains the Feigin-Frenkel center $\fz(\fg)\subset\VV(\fg)_{-h^\vee}$. Denote by  $\VV(\fg)_{\chi=0}$ the quotient of
$\VV(\fg)_{-h^\vee}$ by the ideal generated by $\fz(\fg)$. There is an embedding
\begin{equation}
\label{emd-of-aff}
\VV(\fg)_{\chi=0}\hookrightarrow\Gamma(X,\cD^{ch}_{G/B}),
\end{equation}
which gives the cohomology $H^\bullet(G/B,\cD^{ch}_{G/B})$ a $\VV(\fg)_{\chi=0}-$, hence $\widehat{\fg}-$, module structure. The latter was computed in \cite{AM} to the effect that
\begin{equation}
\label{mod-str-am}
H^i(G/B,\cD^{ch}_{G/B})=\bigoplus_{w\in W,\,l(w)=i}\VV(\fg)_{\chi=0},
\end{equation}
where $W$ is the Weyl group, $l(.)$ is the length function.
Our task is to compute the vertex algebra structure on $H^\bullet(G/B,\cD^{ch}_{G/B})$, see
Lemma~\ref{vert-alg-str-on-coho}.

Denote by $\text{Sing}$ the commutator (i.e., the space of singular vectors) of
$\VV(\fg)_{\chi=0}$ in $H^\bullet(G/B,\cD^{ch}_{G/B})$, cf. (\ref{emd-of-aff}).
More explicitly,
\begin{equation}
\label{def-of-sing}
\text{Sing}=\{a\in H^\bullet(G/B,\cD^{ch}_{G/B})\text{ s.t. }g_{(n)}a=0\, \forall
g\in \VV(\fg)_{\chi=0},n\geq 0\}.
\end{equation}
It is an easy consequence of (\ref{Borcherds-identity}) that $\text{Sing}$ is a vertex
subalgebra of $H^\bullet(G/B,\cD^{ch}_{G/B})$. Furthermore, the irreducibility of
$\VV(\fg)_{\chi=0}$, \cite{FG3}, implies that $\text{Sing}$ is a $\CC$-span of several
copies of the vacua, one for each copy of $\VV(\fg)_{\chi=0}$ appearing in the decomposition (\ref{mod-str-am}). It follows from decomposition (\ref{mod-str-am}) that
the two embeddings $\text{Sing}\hookrightarrow H^\bullet(G/B,\cD^{ch}_{G/B})
\hookleftarrow \VV(\fg)_{\chi=0}$ give
a vertex algebra isomorphism
\begin{equation}
\label{main-res-form-0}
\text{Sing}\otimes_{\CC} \VV(\fg)_{\chi=0}\iso H^\bullet(G/B,\cD^{ch}_{G/B}).
\end{equation}
To determine a vertex algebra structure on $H^\bullet(G/B,\cD^{ch}_{G/B})$, we need that on $\text{Sing}$. Now let us recall that a vertex algebra $V$ is called {\it commutative} if multiplications $_{(n)}=0$ if $n\geq 0$. In this case, the triple $(V,_{(-1)},\vacu)$ is a unital associative commutative algebra with derivation
\begin{equation}
\label{def-of-deriv}
T:V\rightarrow V,\;T(a)=a_{(-2)}\vacu.
\end{equation}
 In fact, the assignment $V\mapsto (V,_{(-1)},\vacu, T)$ establishes
an equivalence between the category of commutative vertex algebras and associative commutative unital algebras with derivation, see e.g. \cite{K}. In particular, any choice of a derivation makes a commutative associative unital algebra into a vertex algebra.

An associative commutative algebra we wish to focus on is $H^\bullet(\fn,\CC)$, the cohomology of the maximal nilpotent $\fn\subset\fg$ with trivial coefficients.

\begin{thm}
\label{main-thm}

(i)   The vertex algebra $\text{Sing}$ is commutative and, for some  derivation of $H^\bullet(\fn,\CC)$, there is a vertex algebra isomorphism
\begin{equation}
\label{sing-coho}
\text{Sing}\iso H^\bullet(\fn,\CC).
\end{equation}

(ii) $H^\bullet(\fn,\CC)$ being equipped with a vertex algebra stucture as in (\ref{sing-coho}), there is
a vertex algebra isomorphism
\begin{equation}
\label{main-res-form}
H^\bullet(\fn,\CC)\otimes_{\CC} \VV(\fg)_{\chi=0}\iso H^\bullet(G/B,\cD^{ch}_{G/B}).
\end{equation}

(iii) If $\fg$ is either $sl_2$ or $sl_3$, then the derivation of $ H^\bullet(\fn,\CC)$
chosen in assertion (i) is 0.
\end{thm}

\begin{cor}
\label{ass-on-centre}
The center of $H^\bullet(G/B,\cD^{ch}_{G/B})$ is isomorphic to $H^\bullet(\fn,\CC)$ as an
associative commutative algebra.
\end{cor}

{\em Proof of Corollary~\ref{ass-on-centre}.} Recall that the center of a vertex algebra
$V$ is $\{v\in V:\; v_{(n)}V=0\,\forall n\geq 0\}$. We have $\text{Sing}_{(n)}\VV(\fg)_{\chi=0}$, $n\geq 0$ by definition, $\text{Sing}_{(n)}\text{Sing}=0$ by virtue of Theorem~\ref{main-thm} (i), hence
$\text{Sing}_{(n)}(\text{Sing}\otimes\VV(\fg)_{\chi=0})=0$; it remains to use (\ref{sing-coho}) and (\ref{main-res-form}). $\Box$

\bigskip

{\em Proof of Theorem~\ref{main-thm}.}

{\em Proof of assertion (i).}

{\bf A.} Let us establish a semi-infinite cohomology interpretation of $\text{Sing}$; this will serve as a bridge
to the Lie algebra cohomology $H^\bullet(\fn,\CC)$.
Denote by $\cC^{\infty/2+\bullet}(L\fn,\CC)$ the semi-infinite cochain complex for the loop
algebra $L\fn$, i.e., the vertex algebra whose various differentials lead to various cohomology
theories of $L\fn$. For $V$  a vertex algebra and $V(\fn)\rightarrow V$ a vertex algebra morphism,
$\chi:\fn\rightarrow\CC$ the principal character, let
\[
\cC^{\infty/2+\bullet}(L\fn; V)_{DS}\stackrel{\text{def}}{=}(\cC^{\infty/2+\bullet}(L\fn,\CC)\otimes V, d_{DS})
\]
 be the Drinfeld-Sokolov reduction complex. It is a differential graded vertex algebra with differential $d_{DS}$;
 its cohomology will be denoted by $H^{\infty/2+\bullet}(L\fn; V)_{DS}$.

 Suppose now $V=H^\bullet(G/B,\cD^{ch}_{G/B})$. A singular vector $v\in\text{Sing}$ determines a cohomology class:
\begin{equation}
\label{cocycles}
\text{Sing}\rightarrow  H^{\infty/2+0}(L\fn; H^\bullet(G/B,\cD^{ch}_{G/B}))_{DS},\; v\mapsto\text{ class of }\vacu\otimes v.
\end{equation}
Frenkel and Gaitsgory \cite{FG3} have computed the cohomology $H^{\infty/2+\bullet}(L\fn; \VV(\fg)_{\chi=0})_{DS}$. Since
$H^\bullet(G/B,\cD^{ch}_{G/B}))_{DS}$ is a direct sum of several copies of the   module $\VV(\fg)_{\chi=0}$, (\ref{mod-str-am}),
their result applies and implies that map (\ref{cocycles}) defines a vertex algebra isomorphism
\begin{equation}
\label{sin-via-semi-inf}
\text{Sing}\iso  H^{\infty/2+\bullet}(L\fn; H^\bullet(G/B,\cD^{ch}_{G/B}))_{DS},\; v\mapsto\text{ class of }\vacu\otimes v.
\end{equation}
It is convenient to recast this result as follows. Consider a sheaf of differential graded vertex algebras
$\cC^{\infty/2+\bullet}(L\fn; \cD^{ch}_{G/B})_{DS}$. The {\it hypercohomology} $H^\bullet(G/B,\cC^{\infty/2+\bullet}(L\fn; \cD^{ch}_{G/B})_{DS})$
arises; this is a vertex algebra by virtue of (an obvious version of) Lemma~\ref{vert-alg-str-on-coho}. Isomorphism (\ref{sin-via-semi-inf})
implies
\begin{equation}
\label{sin-via-semi-inf-geometr}
\text{Sing}\iso  H^\bullet(G/B,\cC^{\infty/2+\bullet}(L\fn; \cD^{ch}_{G/B})_{DS}).
\end{equation}
Let us now relate this to the ordinary Lie algebra cohomology $H^\bullet(\fn,\CC)$.

{\bf B.} For an $\fn$-module $V$, denote by
$\cC^\bullet(\fn,V)$ the standard cochain complex with values in $V$. Since $\fn\subset\fg$ operates on $G/B$, a sheaf
of commutative differential graded algebras $\cC^\bullet(\fn,\cO_{G/B})$ arises. We have an obvious
embedding $\cO_{G/B}\hookrightarrow\cD^{ch}_{G/B}$, $f\mapsto f\vacu$. We would like to assert that it is a vertex algebra
morphism, but it is not because $\cO_{G/B}$ does not carry an appropriate derivation because if $f\in\cO_{G/B}\subset\cD^{ch}_{G/B}$, then $f_{(-2)}1\not\in\cO_{G/B}$, cf. (\ref{def-of-deriv}).
We can, however, equip $\cO_{G/B}$ with multiplications $_{(n)}$, $n\geq -1$, by declaring that $_{(-1)}$ is the standard multiplication
on $\cO_{G/B}$ and $_{(n)}=0$ if $n\geq0$. Let us for the purpose of the remainder of this proof call a vector space space with
multiplications $_{(n)}$, $n\geq -1$, an {\em algebra}. Thus the morphism $\cO_{G/B}\hookrightarrow\cD^{ch}_{G/B}$ becomes an {\em algebra}
morphism.

It naturally extends to a differential graded
 algebra morphism
\begin{equation}
\label{usual-semi-inf}
\cC^\bullet(\fn,\cO_{G/B})\rightarrow \cC^{\infty/2+\bullet}(L\fn; \cD^{ch}_{G/B})_{DS}.
\end{equation}
Since $H^\bullet(G/B,\cO_{G/B})=H^0(G/B,\cO_{G/B})=\CC$, $H^\bullet(G/B,\cC^\bullet(\fn,\cO_{G/B}))=H^\bullet(\fn,\CC)$.
Therefore, (\ref{usual-semi-inf}) gives
\begin{equation}
\label{usual-semi-inf-coho-level}
H^\bullet(\fn,\CC)\rightarrow H^\bullet(G/B,\cC^{\infty/2+\bullet}(L\fn; \cD^{ch}_{G/B})_{DS}).
\end{equation}

\begin{lem}
\label{ordin-semi-inf-iso}
Map (\ref{usual-semi-inf-coho-level}) is an algebra isomorphism.
\end{lem}

This lemma implies Theorem~\ref{main-thm} (i)
at once because assertion (i) is the lemma combined with  (\ref{sin-via-semi-inf-geometr}); in particular,
the derivation of $H^\bullet(\fn,\CC)$, the missing ingredient, is one  pulled back from $\text{Sing}$.

{\em Proof of Lemma~\ref{ordin-semi-inf-iso}.} We shall use the approach of \cite{AM}.
 Starting with the BGG-resolution of the trivial $\fg$-module
\begin{eqnarray}
&0\rightarrow\CC\rightarrow M^\bullet\label{bgg}\\
&M^\bullet:\; M_0^c\rightarrow\bigoplus_{w\in W, l(w)=1}M^c_{w\cdot 0}\rightarrow\cdots\rightarrow\bigoplus_{w\in W, l(w)=i}M^c_{w\cdot 0}\rightarrow\cdots,
\nonumber
\end{eqnarray}
we pull the bootstraps to obtain, first, the Cousin-Grothendieck resolution of $\cO_{G/B}$, via the Beilinson-Bernstein localization functor $\Delta$,
\begin{equation}
\label{cous-groth}
0\rightarrow\cO_{G/B}\rightarrow \Delta(M^\bullet);
\end{equation}
second, a resolution of $\cD^{ch}_{G/B}$, using a version of the Zhu functor
\cite{AChM}
\begin{equation}
\label{cous-groth-chir}
0\rightarrow\cD^{ch}_{G/B}\rightarrow \cZ hu\circ\Delta(M^\bullet).
\end{equation}
In the derived category, the map $\cO_X\rightarrow\cD^{ch}_{G/B}$ is nothing but the natural embedding
\[
\Delta(M^\bullet)\hookrightarrow  \cZ hu\circ\Delta(M^\bullet),
\]
hence
\begin{equation}
\label{reformul}
\cC^\bullet(\fn,\Delta(M^\bullet))\rightarrow  \cC^{\infty/2+\bullet}(L\fn;\cZ hu\circ\Delta(M^\bullet))_{DS}.
\end{equation}
The map of hypercohomology induced by the latter is precisely (\ref{usual-semi-inf-coho-level}), but now we  can compute
it in a different way.
An application of $R\Gamma(G/B,.)$ to both sides of (\ref{reformul}) gives
\[
\cC^\bullet(\fn,M^\bullet)\rightarrow  \cC^{\infty/2+\bullet}(L\fn;\Gamma(G/B,\cZ hu\circ\Delta(M^\bullet)))_{DS}.
\]
Since $M^\bullet$ is a co-free resolution of $\CC$, the cohomology of the left hand side is $H^\bullet(\fn,\CC)$.
Furthermore, the cohomology classes are 0-cochains defined by highest weight vectors of the contragredient Verma
modules $M^c_{w\cdot 0}$, cf. (\ref{bgg}).

Similarly, it is the essence of the proof of (\ref{mod-str-am}) proposed in \cite{AM} that the cohomology
of the right hand side is $H^\bullet(G/B,\cC^{\infty/2+\bullet}(L\fn; \cD^{ch}_{G/B})_{DS})$ and that the cohomology
classes are precisely 0-cochains defined by the highest weight vectors of $\Gamma(G/B,\cZ hu\circ\Delta(M^\bullet))$.

The map $M^\bullet\rightarrow \Gamma(G/B,\cZ hu\circ\Delta(M^\bullet))$, by its definition, sends highest weight vectors
to highest weight vectors. Hence the resulting map
\[
H^\bullet(\fn,\CC)\rightarrow H^\bullet(G/B,\cC^{\infty/2+\bullet}(L\fn; \cD^{ch}_{G/B})_{DS})
\]
is an isomorphism.

{\em Assertion (ii)} is  decomposition (\ref{main-res-form-0}) and assertion (i) combined.

{\em Proof of assertion (iii).} As follows from (\ref{def-of-deriv}),
we need to show that $a_{(-2)}\vacu=0$ for $\forall\,a\in\text{Sing}$. Definition (\ref{def-of-sing}) implies
that both $a$ and $a_{(-2)}\vacu$ are highest weight vectors w.r.t. $\widehat{\fg}$ of conformal weights that differ from each other by 1.
 On the other hand, if $a\in H^{i}(G/B,\cD^{ch}_{G/B})$,
then $a_{(-2)}\vacu\in H^{i}(G/B,\cD^{ch}_{G/B})$ too. However, Theorem 1.1 (2) of \cite{AM} shows at once that if
 $\fg$ is either $sl_2$ or $sl_3$, then all highest weight vectors of the same cohomological degree have the same conformal weight; hence
 $a_{(-2)}\vacu$ must be 0.

This concludes the proof of  Theorem~\ref{main-thm}. $\Box$

\section{Bosonization over the projective line}
\label{Bosonization over the projective line}
\subsection{The CDO}
\label{the cdo}
The CDO $\cD^{ch}_{\prline}$ can be defined very explicitly.
Let $\{\CC_0,\CC_\infty\}$ be the standard atlas of $\prline$
with $x$, $y=1/x$ the respective coordinates. The space
$\cD^{ch}_{\prline}(\CC_0)$ is what is known as {\it the $\beta\gamma$-system},
a vertex algebra generated by 2 even elements, $x$, $\partial_x$, and relations
\begin{equation}
\label{def-b-g}
\partial_{x (n)}x=-x_{(n)}\partial_x=\delta_{n0}\vacu\text{ if }n\geq 0.
\end{equation}
$\cD^{ch}_{\prline}(\CC_\infty)$ is defined similarly with $x$ replaced with $y$.
That $x$ and $y$ play a 2-faceted role, a coordinate and a vertex algebra element, should not lead to too much confusion for their transformation laws, which we are about to introduce,
are the same in each of their capacities.

One can localize $\cD^{ch}_{\prline}(\CC_0)$ ($\cD^{ch}_{\prline}(\CC_\infty)$)
over $\CC_0$  ($\CC_\infty$ resp.) and then glue over $\CC^*=\CC_0\cap\CC_\infty$
sending
\begin{equation}
\label{transf-form}
x\mapsto 1/y,\partial_x\mapsto -\partial_{y (-1)}(y_{(-1)}y)-2y_{(-2)}\vacu.
\end{equation}
The result is the $\cD^{ch}_{\prline}$.

\subsection{The Friedan-Martinec-Shenker bosonization}
\label{the bosonization}
Consider the lattice $L=\ZZ p\oplus\ZZ q$ with $(p,p)=-(q,q)=1$, $(p,q)=0$.
There arises a {\it lattice vertex algebra}, $V_L$, see an exposition in\cite{K}, sect.5.4,5.5. Its main ingredients are the Lie alhebra $\hat{\fh}$ with basis
$\{p_n,q_n,\;n\in\ZZ\}$ and commutation relations
\[
[p_n,p_m]=-[q_n,q_m]=n\delta_{n,-m},\;[p_n,q_m]=0,
\]
$S$, the vacuum $\hat{\fh}$-module, and $\CC[L]_{\epsilon}$, a
twisted group algebra of $L$ with multiplication
\[
e^\alpha e^\beta=\epsilon(\alpha,\beta)e^{\alpha+\beta},
\]
see \cite{K}, sect.5.5, for an explicit formula for the cocycle $\epsilon$.

As a vector space, $V_L=\CC[L]_{\epsilon}\otimes S$. The formulas for products
$_{(n)}$ are classic and can be found in $\it loc. cit.$; we shall only record
the generating function for the operators $e^\alpha_{(n)}$ given by the celebrated
{\it vertex operators}:
\[
\sum_{n\in\ZZ}e^\alpha_{(n)}z^{-n-1}=e^\alpha z^{(\alpha,.)}\exp{
\sum_{n<0}\frac{\alpha_n}{-n}z^{-n}}\exp{
\sum_{n>0}\frac{\alpha_n}{-n}z^{-n}},
\]
where $(\alpha,.)$ is understood as a function that sends $e^\beta\otimes s$ to
$(\alpha,\beta)e^\beta\otimes s$.

Define the grading $V_L=\oplus_{m\in\ZZ}V_L[m]$, where $V_L[m]$ is the eigenspace
of $(p+q,.)$ with eigenvalue $m$. We let
\[
V_{L\geq}=\oplus_{m\geq 0}V_L[m],\;V_{L\leq}=\oplus_{m\leq 0}V_L[m].
\]
If we let $D_\geq=e^p_{(0)}$, $D_\leq=e^{-p}_{(0)}$, then we obtain 2 differential
graded vertex algebras
\[
(V_{L\geq}, D_\geq),\;(V_{L\leq}, D_\leq)
\]
and the cohomology vertex algebras, $H^\bullet(V_{L\geq})$ and $H^\bullet(V_{L\leq})$.

Now consider $\cD^{ch}_{\prline}(\CC_0)$ as a differential graded vertex algebra
with trivial grading and zero differential.
\begin{lem}
\label{fr-m-sh-f-f}
The assignment
\[
x\mapsto e^{p+q},\;\partial_x\mapsto p_{(-1)}e^{-p-q}
\]
defines a quasiisomorphism of differential graded vertex algebras
\[
\rho_\geq:\;\cD^{ch}_{\prline}(\CC_0)\rightarrow (V_{L\geq}, D_\geq).
\]
In particular,
\[
\cD^{ch}_{\prline}(\CC_0)\stackrel{\sim}{\rightarrow}H^0(V_{L\geq}).
\]
\end{lem}

This is known as Friedan-Martinec-Shenker bosonization, \cite{FMS}. The assertion
about $H^0(V_{L\geq})$ was proved in \cite{FF}. The vanishing of the higher cohomology
(undoubtedly known to Feigin and Frenkel) is  easy and follows from the fact that
$e^{-p}_{(-1)}$ is a contracting homotopy; indeed $[e^p_{(0)},e^{-p}_{(-1)}]=1$.

The advantage of passing from $\cD^{ch}_{\prline}(\CC_0)$ to $(V_{L\geq}, D_\geq)$ is that
in the latter the element $x$ is invertible: just let $x^{-1}=e^{-p-q}$. Upon localization,
we obtain a vertex  algebra morphism
\[
\cD^{ch}_{\prline}(\CC^*)\rightarrow V_L[0],
\]
 which is readily seen
to be an isomorphism. Therefore, an embedding
$\cD^{ch}(\CC_\infty)\hookrightarrow V_L[0]$, compatible with
$\cD^{ch}(\CC_0)\hookrightarrow V_L[0]$ that appears in Lemma~\ref{fr-m-sh-f-f}, must arise.
It does and its construction is natural: the map $x\mapsto x^{-1}$ lifts to an automorphism of $\cD^{ch}(\CC^*)$. In terms of $V_L[0]$, the latter is nothing but the automorphism
engendered by the lattice automorphism $p\mapsto-p$, $q\mapsto -q$. Hence the following assertion, a companion of Lemma~\ref{fr-m-sh-f-f}.
\begin{cor}
\label{cor-of-fr-m-sh-f-f}
The assignment
\[
y\mapsto e^{-p-q},\;\partial_y\mapsto -p_{(-1)}e^{p+q}
\]
defines a quasiisomorphism of differential graded vertex algebras
\[
\rho_\leq:\;\cD^{ch}_{\prline}(\CC_\infty)\rightarrow (V_{L\leq}, D_\leq)
\]
so that
\[
\rho_\geq|_{\cD^{ch}_{\prline}(\CC_0)\cap \cD^{ch}_{\prline}(\CC_\infty)}=
\rho_\leq|_{\cD^{ch}_{\prline}(\CC_0)\cap \cD^{ch}_{\prline}(\CC_\infty)},
\]
where the intersection $\cD^{ch}_{\prline}(\CC_0)\cap \cD^{ch}_{\prline}(\CC_\infty)$ is that
of the images w.r.t. the embeddings $\cD^{ch}_{\prline}(\CC_0)\hookrightarrow \cD^{ch}_{\prline}(\CC^*)\hookleftarrow
\cD^{ch}_{\prline}(\CC_\infty)$
In particular,
\[
\cD^{ch}_{\prline}(\CC_\infty)\stackrel{\sim}{\rightarrow}H^0(V_{L\leq}).
\]
\end{cor}

\subsection{A bosonization \`a la Borisov}
\label{The bosonozation a la Borisov}
Consider a complex
\[
0\rightarrow V_L[0]\stackrel{(D_\geq,D_\leq)}{\longrightarrow}V_L[1]\oplus V_L[-1]\rightarrow 0.
\]
It follows at once from Lemma~\ref{fr-m-sh-f-f} and Corollary~\ref{cor-of-fr-m-sh-f-f} that its 0-th cohomology
is $H^0(\prline,\cD^{ch}_{\prline})$. We shall now extend it to a complex that computes the entire $H^\bullet(\prline,\cD^{ch}_{\prline})$.
Our argument is a straightforward purely even version of the construction that Borisov proposed on the case of the chiral
de Rham complex \cite{B}.

Consider the \u Cech complex  $\check{\cC}^\bullet=\check{\cC}^\bullet(\fU,\cD^{ch}_{\prline})$
 attached to the open cover $\fU=\{\CC_0,\CC_\infty\}$:

\begin{equation}
\xymatrix{ \cD^{ch}_{\prline}(\CC^*)\\
           \cD^{ch}_{\prline}(\CC_0)\oplus\cD^{ch}_{\prline}(\CC_\infty)\ar[u]}
\end{equation}
Lemma~\ref{fr-m-sh-f-f} and Corollary~\ref{cor-of-fr-m-sh-f-f} give us a resolution of $\check{\cC}^\bullet$, to be denoted
$\cR^\bullet(\check{\cC}^\bullet)$, that is a double complex as follows:
\begin{equation}
\label{bor-res-double-compl}
\xymatrix{ V_L[0]\ar[r]&0\ar[r]&0\ar[r]&\cdots\\
           V_L[0]\oplus V_L[0]\ar[u]\ar[r]^<<<<<{(D_\geq,D_\leq)}&V_L[+1]\oplus V_L[-1]\ar[u]\ar[r]^{(D_\geq,D_\leq)}&V_L[+2]\oplus V_L[-2]
           \ar[u]\ar[r]
           &\cdots
           }
\end{equation}

Recall that the cup-product formula (\ref{mult-cech-coch}) makes $\check{\cC}^\bullet$ into a differential graded algebra with
products $_{(n)}$, $n\in\ZZ$, and it is this algebra structure that defines a vertex algebra structure on $H^\bullet(\prline,\cD^{ch}_{\prline})$,
Lemma~\ref{vert-alg-str-on-coho}. The cup-product is easy to extend from $\check{\cC}^\bullet$  to $\cR^\bullet(\check{\cC}^\bullet)$:
in order to do that simply write $\cR^\bullet(\check{\cC}^0)_i=\cR^\bullet(\check{\cC}^0_i)$, $\cR^\bullet(\check{\cC}^1)_{01}=\cR^\bullet(\check{\cC}^0_{01})$ and then use the same (\ref{mult-cech-coch}). The quasiisomorphism
\[
\check{\cC}^\bullet\rightarrow \cR^\bullet(\check{\cC}^\bullet)
\]
becomes a differential graded algebra quasiisomorphism, the quasiisomorphism assertion being the result of computing the cohomology
of $\cR^\bullet(\check{\cC}^\bullet)$  w.r.t. the horizontal differential in (\ref{bor-res-double-compl}).

A computation starting with the vertical differential will again display a collapse of the spectral sequence and show that
$H^\bullet(\prline,\cD^{ch}_{\prline})$ equals the cohomology of the complex
\begin{equation}
\label{bor-res-single-compl}
\xymatrix{
           V_L[0]\ar[r]^<<<<<{(D_\geq,D_\leq)}&V_L[+1]\oplus V_L[-1]\ar[r]^{(D_\geq,D_\leq)}&V_L[+2]\oplus V_L[-2]
           \ar[r]
           &\cdots
           }
\end{equation}
This complex is also a differential graded algebra with a family of products $_{(n)\tilde{ }}$. It is a bit unexpected
that not only this algebra structure defines a vertex algebra structure on the cohomology, but it itself is a vertex algebra
on the nose. To see this, note that the cup-product on 0-chains being defined by $(f_{(n)\tilde{ }}g)(i)=f(i)_{(n)\tilde{ }}g(i)$,
cf. (\ref{mult-cech-coch}), the products on (\ref{bor-res-single-compl}) are
\begin{equation}
\label{degen-prod}
u_{(n)\tilde{ }}v=\left\{\begin{array}{cc}u_{(n)}v&\text{ if }u\in V_L[m],v\in V_L[n], m\cdot n\geq 0\\
0&\text{ if }u\in V_L[m],v\in V_L[n], m\cdot n< 0
\end{array}\right.
\end{equation}
 It is easy to show directly that (\ref{degen-prod}) is a vertex algebra structure; it is
even more useful to obtain it as a degeneration of one on $V_L$. Indeed upon rescaling
\[
V_L[m]\ni v\mapsto t^{|m|}v,
\]
the product $_{(n)}$ is replaced with $_{(n),t}$ so that
\[
u_{(n),t}v=\left\{\begin{array}{cc}
u_{(n)}v &\text{ if }\text{deg}u\cdot\text{deg}v\geq 0\\
t^{|\text{deg}u|+|\text{deg}v|-|\text{deg}u+\text{deg}v|}u_{(n)}v &\text{ if }\text{deg}u\cdot\text{deg}v< 0
\end{array}\right.
\]
This gives a family of vertex algebras, $\{V_L^t,t\in\CC^*\}$, isomorphic to each other
 and with a well-defined limit as $t\rightarrow 0$. Furthermore,
$V_L^0\stackrel{\text{def}}{=}\lim_{t\rightarrow 0}V_L^t$ is precisely (\ref{degen-prod});
hence the latter  also defines a vertex algebra structure. Since $D_\geq+D_\leq$ is a derivation of the latter, our discussion proves the following.
\begin{lem}
There is a vertex algebra isomorphism
\label{coho-a la-bor-alg}
\[
H^\bullet(\prline,\cD^{ch}_{\prline})\stackrel{\sim}{\rightarrow}
H^{\bullet}_{D_\geq+D_\leq}(V_L^0).
\]
\end{lem}
Since $D_\geq+D_\leq$ remains a derivation in the deformed algebra, we can consider
$H^{\bullet}_{D_\geq+D_\leq}(V_L^t)$.
\begin{thm}
\label{deformed-coho}
\[
H^{\bullet}_{D_\geq+D_\leq}(V_L^t)=\left\{
\begin{array}{cc}
H^\bullet(\prline,\cD^{ch}_{\prline})&\text{ if }t=0\\
0&\text{ otherwise.}
\end{array}
\right.
\]
\end{thm}

{\em Proof.} The argument is similar to that in \cite{MS}, and we shall be brief.
The differential on $V_L^t$ can be written as $D=D_0+tD_1$, where $D_0$ coincides
with the differential on the degenerated $V_L^0$, $D_1$ commutes with $D_0$ and maps
in the opposite direction. Indeed, while $D_0$ is as in (\ref{bor-res-single-compl}),
the complex $(V_L^t,D_0+tD_1)$ is as follows:
\begin{equation}
\label{bor-res-single-compl-deformed}
\xymatrix{
           V_L[0]\ar@<1ex>[r]^<<<<<{(D_\geq,D_\leq)}&V_L[+1]\oplus V_L[-1]\ar@<1ex>[r]^{(D_\geq,D_\leq)}\ar@<1ex>[l]^>>>>>>{t(D_\leq,D_\geq)}& V_L[+2]\oplus V_L[-2]
           \ar@<1ex>[r]\ar@<1ex>[l]^{t(D_\leq,D_\geq)}
           &\cdots
           \ar@<1ex>[l]
           }
\end{equation}
The case $t=0$ of the theorem is but Lemma~\ref{coho-a la-bor-alg}. If $t\neq 0$, then
a spectral sequence arises with $E_1$ the $D_0$-cohomology. The latter is
known, \cite{MSV,AM}: $E_1=H^0(\prline,\cD^{ch}_{\prline})\oplus H^1(\prline,\cD^{ch}_{\prline})$, $H^0(\prline,\cD^{ch}_{\prline})= H^1(\prline,\cD^{ch}_{\prline})=\VV(sl_2)_{\chi=0}$, cf. (\ref{mod-str-am}).
Therefore, $(E_1, d_1)$ is
\begin{equation}
\label{e-1}
 tD_1: H^1(\prline,\cD^{ch}_{\prline})\rightarrow H^0(\prline,\cD^{ch}_{\prline}).
\end{equation}
 $tD_1$ is a $\hat{\fg}$-module morphism.
Let us see how it acts on the highest weight vector.

The highest weight vector of $H^1(\prline,\cD^{ch}_{\prline})$ was computed in \cite{MSV};
it equals  the \u Cech cochain $\CC^*\mapsto x_{(-2)}x^{-1}$. According to Lemma~\ref{fr-m-sh-f-f}, upon bosonization it is identified with
\[
e^{-p-q}_{(-1)}((p+q)_{(-1)}e^{p+q})=(p+q)_{(-1)}\vacu\in V_L[0].
\]
This places our cocycle in the uppermost-leftmost corner of (\ref{bor-res-double-compl}).
In order to identify it with a cocycle in $V_L^0$, see (\ref{bor-res-single-compl}), we have to do a bit of a diagram chase.

It is clear that the cocycle belongs to the image of $((p+q)_{(-1)}\vacu,0)$ under the action of the vertical differential in (\ref{bor-res-double-compl}), hence cohomologous to
\[
(D_\geq((p+q)_{(-1)}\vacu),0)=(e^p_{(0)}((p+q)_{(-1)}\vacu),0)=(-e^p,0)\in V_L^0[1].
\]
Therefore, its image under the action of $tD_1$ is, see (\ref{bor-res-single-compl-deformed}),
\[
tD_1((-e^p,0))=t(-D_\leq e^p,0)=te^{-p}_{(0)}e^p=t\vacu,
\]
a highest weight vector of $H^0(\prline,\cD^{ch}_{\prline})$ if $t\neq 0$. Since
$H^0(\prline,\cD^{ch}_{\prline})$ and $H^1(\prline,\cD^{ch}_{\prline})$ are irreducible,
morphism (\ref{e-1}) is an isomorphism. Therefore, $E_2=0$, and the assertion of Theorem~\ref{deformed-coho} follows. $\Box$

\begin{rem}
\label{conn-to-witt}
Witten explains \cite{W} that $H^\bullet(\prline,\cD^{ch}_{\prline})=H^\bullet_{D_\geq+D_\leq}(V_L^0)$ is
Witten's half-twisted algebra of a certain $\sigma$-model with (0,2)-supersymmetry in perturbative regime.
If we are to believe that  deforming $H^\bullet_{D_\geq+D_\leq}(V_L^0)$ to $H^\bullet_{D_\geq+D_\leq}(V_L^t)$
means taking into account instanton corrections -- such apparently was original Borisov's intention \cite{B}, and indeed this
construction gives quantum cohomology in case of $\PP^n$ and the chiral de Rham complex \cite{MS}, see also \cite{FL} -- then the vanishing
$H^\bullet_{D_\geq+D_\leq}(V_L^t)=0$ if $t\neq 0$ confirms Witten's assertion \cite{W} that nonperturbatively half-twisted
algebra vanishes.
\end{rem}

%


\section*{Acknowledgments}

FM participated in the 8th International Workshop
"Lie Theory and Its Applications in Physics" held from
15--21 June 2009 in  Varna, Bulgaria and would like to thank the organizers
for exceptional hospitality.


\footnotesize{T.A: Department of Mathematics, Nara Women's
University, Nara 630-8506, Japan. E-mail address:
arakawa@cc.nara-wu.ac.jp

F.M.: Department of Mathematics, University of Southern California, Los Angeles, CA 90089, USA.
E-mail:fmalikov@usc.edu}
\end{document}